\documentclass[11pt]{amsart}

\usepackage{amsmath}
\usepackage{amsfonts}
\usepackage{amssymb}
\usepackage{amsthm}
\usepackage{graphics}
\usepackage{amscd}
\usepackage{graphicx,epsfig}
\usepackage{color}
\usepackage[all]{xy}
\usepackage{url}

\renewcommand{\epsilon}{\varepsilon}

\newcommand{\boL}{\mathcal{L}}

\newcommand{\boN}{\mathcal{N}}

\newcommand{\R}{\mathbb{R}}

\newcommand{\N}{\mathbb{N}}
\newcommand{\la}{\langle}
\newcommand{\ra}{\rangle}
\newcommand{\eps}{\varepsilon}

\newtheorem{prop}{Proposition}

\renewcommand{\phi}{\varphi}

\newtheorem*{thm*}{Theorem}

\newcounter{remark}

\newcounter{case}

\newcounter{construction}

\newcounter{fact}

\newcommand{\T}{\mathbb{T}}

\title{Erratum to "Minimal surfaces in finite volume noncompact hyperbolic $3$-manifolds"}
\author{Pascal Collin, Laurent Hauswirth, Laurent Mazet, Harold Rosenberg}

\begin{document}

\maketitle

In this note, we explain how a mistake made in \cite{CoHaMaRo} can be corrected. Actually 
this mistake appears in the proof of Proposition~8 (the second maximum principle) and 
was brought to our attention by A.~Song~\cite{Son2}. Let us notice that unfortunately, 
we did not find an alternative proof of this proposition but we found an alternative 
proposition. The new proposition does not change the subsequent applications we made 
of the original proposition. At the end of the note we explain which modifications should 
be done where 
the original Proposition~8 is applied. 

The difference with the original Proposition~8 is that here we have to assume a control on 
the index of the minimal surface $\Sigma$.

\subsection*{The new proposition}

Let  $(\T,ds_\T^2)$ be a flat $2$ torus of diameter $1$ and consider 
$M=\T\times 
[-\frac12,+\infty)$ endowed with the metric $g=e^{-2h(t)}\Lambda^2 
ds_{\T}^2+dt^2$ ; here $\Lambda$ is a positive constant and $h$ satisfies the 
following assumptions:
\begin{itemize}
\item $h(0)=0$
\item $h'(t)>0$
\item $h'$ is bounded in the $C^k$ norm.
\end{itemize}

The main example is $h(t)=t$. We denote by $T_s=\T\times\{s\}$.

\begin{prop}\label{prop}
Let $i_0\in\N$. There is a $\Lambda_0$ such that the following is true. For any 
$\Lambda\le \Lambda_0$ 
and any compact embedded minimal surface $\Sigma\subset M$ with $\partial 
\Sigma\subset T_{-\frac12}$ and index less than $i_0$, we have $\Sigma\subset 
\T\times[-\frac12,0]$
\end{prop}

\begin{proof}
If the proposition is not true there is a sequence $\Lambda_n\to 0$ and 
$\Sigma_n$ a minimal surface in $M$ with index less than $i_0$, boundary in 
$T_{-\frac12}$ and the maximum of the function $t$ on $\Sigma_n$ is $t_n>0$.

In order to study this sequence, we translate in the $t$ direction by $-t_n$ 
and make a homothety by $\lambda_n= \frac{e^{h(t_n)}}{\Lambda_n}\to +\infty$. We 
thus 
obtain a compact minimal surface $S_n$ in $\T\times 
[-(\frac12+t_n)\lambda_n ,0]$ 
endowed with the metric $g_n=e^{-2h_n(t)}
ds_{\T}^2+dt^2$ (where $h_n(t)=h(\frac t{\lambda_n}+t_n)-h(t_n)$) with boundary 
in $T_{-(\frac12 +t_n)\lambda_n}$ and containing 
a point in $T_0$. We notice that $h_n\to 0$ uniformly on any compact. So the ambient 
space smoothly converges to 
$X=\T\times \R_-$ endowed with the flat metric $ds_{\T}^2\times dt^2$.

Let $\epsilon$ be positive. Since $S_n$ has index at most $i_0$, there is a set 
$E^n(\eps)$ 
of at most $i_0$ 
points 
in $S_n$ such that, for any $p\in 
X$ with $d(p,E^n(\eps))>2\epsilon$, $S_n\cap B(p,\epsilon)$ is stable. Actually, by the 
work of Chodosh, Ketover and Maximo (Theorems~1.5 and 1.17 in \cite{ChKeMa}), up to a 
subsequence, there is a set $E^n$ of at 
most $i_0$ points in $X$ and a constant $\kappa > 0$ such that for all $x\in S_n$,
\[
| A_n|(x) \min\{ 1, d(x, E^n)\} < \kappa,
\]
where $A_n$ is the second fundamental form of $S_n$, $d$ is distance in $X$.

Passing to a subsequence,  the sets $E^n$ converge to a set $E^\infty =  \{\bar 
p_1,\dots,\bar p_j\}$ in $X$ ($j\le i_0$).
On the complement of $E^\infty$, a subsequence of the surfaces $S_n$ converge 
to a minimal 
lamination $\boL$ of $X\setminus E^\infty$.  Moreover $\boL$ extends smoothly to 
$E^\infty$; in any compact ball of $X$, the leaves of $\boL$ have bounded curvature.

%
%

Let us write $\bar p_j=(\bar q_j,\bar t_j)$. Let $\eps_0>0$ be small enough 
such that either $\bar t_j=0$ or $\bar t_j<-\eps_0$. In the following, we will 
consider that the first possibility is true for any $j$ since we won't look at 
the lamination close to the other points.

We notice that, since 
$S_n\cap T_0\neq \emptyset$, we have $\boL\cap T_0\neq \emptyset$. Since $T_0$ 
is minimal and $\boL\subset X$, $T_0$ is a leaf of $\boL$. Besides 
the lamination structure of $\boL$ implies that
\[
\forall \delta>0, \exists \eps>0, \text{ such that }|A_\boL|\le \delta\text{ and 
}|\la\nu_\boL,\partial_t\ra|\ge 1-\delta\text{ on }X_\eps
\]
where $\nu_\boL$ is the unit normal to $\boL$ and $X_\eps=\T\times[-\eps,0]$.

So the lamination structure implies that $\eps_0$ can be chosen such that in 
$X_{\eps_0}$, $\boL$ is transverse to $\partial_t$ and in fact almost orthogonal. 
Finally we also have that, for 
any $L>0$, there is $\eps_L>0$ such that any path $\gamma$ in $\boL$ of length less 
than 
$L$ and starting from $X_{\eps_L}$ stays in $X_{\eps_0-\eps_L}$.


For $\eta,\eps>0$, we define $C_j(\eta,\eps)=D(\bar q_j,\eta)\times[-\eps,0]$ ($D(\bar 
q,\eta)$, a geodesic disk in $\T$). Once 
$\eta$ and $\eps$ are fixed we know that we have convergence $S_n\to \boL$ outside 
$\cup_j C_j(\eta,\eps)$. This implies that for large $n$, in $X\setminus \cup_j 
C_j(\eta,\eps)$, $S_n$ is almost orthogonal to $\partial_t$, its second 
fundamental form has small norm and any path $\gamma\subset S_n\cap 
(X_\eps\setminus \cup_j 
C_j(\eta,\eps)$ of length less than $L$ starting from $X_{\eps_{2L}}$ stays in 
$X_{\eps_0-\eps_{2L}}$.

Let us consider a set of generator of $\pi_1(\T\setminus\{\bar q_1,\dots,\bar q_{\bar 
j}\})$ 
and $L$ an upper bound for the length of these generators. Let us fix $\eps\le \min 
(\eps_0,\eps_{10 L})$. Because of the transversality of  $S_n$ with $\partial_t$ and the 
control of its curvature, we can choose $\eta>0$ such that  any connected component 
of $S_n\cap (D(\bar q_j,\eta)\times[-\eps_0,0])$ intersecting $D(\bar 
q_j,\eta)\times[-\eps_0+\eps,-\eps]$ is a topological disk with boundary a loop on 
$\partial D(\bar q_j,\eta)\times[-\eps_0,0]$ winding around $\bar q_j\times\R$.

Let $\bar p=(\bar q,0)\in T_0\subset \boL$ be a point outside of the 
cylinders $C_j(\eta,\eps)$. Because of the convergence $S_n\to \boL$ outside the 
cylinders, 
for large $n$, there is a point $p_n=(\bar q,\bar t_n)\in S_n$ such that 
$\bar t_n\to 0$. Actually we choose $\bar t_n$ as the largest $t$ such that 
$(\bar q, t)\in S_n$ (we recall that $S_n$ is compact).

Now choose a basis of  generators of $\pi_1(T_0\setminus (\cup 
C_j(\eta,\eps)))$ at $\bar q$ of length less than $2L$. Since, outside the 
$D_j(\eta,\eps)$, 
$S_n$ is transverse to $\partial_t$, such a generator $\gamma:s\in[0,1]\mapsto 
(q(s),0)$ can be lifted vertically to a path 
$\gamma_n$ in $S_n$ starting from $p_n$ for $n$ large (since $L$ is fixed and assuming 
$\bar t_n>-\eps$, we can be 
sure, if $n$ is large, that the vertical lift stays in $X_{\eps_0-\eps}$ where $S_n$ is 
transverse to 
$\partial_t$). We write this 
lift as $\gamma_n:s\mapsto (\bar q(s),t_n(s))$. We have three 
possibilities: $\bar t_n=t_n(0)=t_n(1)$, $\bar t_n=t_n(0)>t_n(1)$ or $\bar 
t_n=t_n(0)<t_n(1)$. By the definition of $\bar t_n$ the second possibility is 
impossible. If the third one occurs, we consider the lift of $-\gamma$, which 
gives $\tilde\gamma_n:t\mapsto (\bar q(1-s),\tau_n(s))$. Since $t_n(1)>t_n(0)=\tau_n(0)$ 
we 
must have $\tau_n(s)<t_n(1-s)$ and $\tau_n(1)<t_n(0)=\bar t_n$ which is 
impossible. So $t_n(0)=t_n(1)$ and $\gamma$ lifts as a loop in $S_n$. This 
implies that the component of $S_n\cap((\T\setminus \cup_j D(\bar 
q_j,\eta))\times[-\eps_0,0])$ containing $p_n$ is a vertical graph over 
$T_0\setminus \cup_j C_j(\eta,\eps))$. It is denoted $\widetilde S_n$. We notice that, as 
$n\to \infty$, this graph 
converges, as graph, to $T_0\setminus \cup_j C_j(\eta,\eps)$. This implies that for $n$ 
large the boundary of the graph is made of curves in $X_{\eps}$. 

Let us look at the case where there is no singular point. In this case, $\widetilde 
S_n$ 
is a 
graph over  the entire $T_0$ so it is a compact connected 
component of $S_n$ with no boundary and such a surface can not exist since the mean 
curvature of each $T_t$ points up.

So now we assume that there are singular points.
Let $p=(q,t)$ be a point in $S_n\cap \partial D(\bar q_j,\eta)\times[-\eps,0]$.  Let 
$\gamma$ be a generator  of $\pi_1((T_0\setminus (\cup 
C_j(\eta,\eps))$ at $p$ of length less than $2L$. As above $\gamma$ can be vertically 
lifted 
to $S_n$ with initial point $p$ to $\gamma_n:s\mapsto (\bar q(s),t_n(s))$. As above we 
have three possibilities: $t=t_n(0)=t_n(1)$, $t=t_n(0)>t_n(1)$ or $t=t_n(0)<t_n(1)$. Let us 
consider the third possibility. Since $-\eps<t_n(0)<t_n(1)$, this implies that we can lift 
$\gamma$ starting from $(q,t_n(1))$ and extend $\gamma_n$ on $[1,2]$. We also have 
$\gamma_n(1)<\gamma_n(2)$. Actually we can extend the definition of $\gamma_n$ to 
$\R_+$. We obtain that $S_n$ contains the sequence of points $(q,t_n(k))_k$. This 
sequence 
converges to a point $(q,\tilde t)\in S_n$ where $S_n$ is transverse to $\partial_t$. This 
gives a contradiction with the fact that $S_n$ is compact and then properly embedded. If 
the second possibility $t_n(0)>t_n(1)$ occurs, we do the same argument as above with 
$-\gamma$. So finally, we have $t_n(0)=t_n(1)$. This implies that $T_0\setminus \cup_j 
C_j(\eta,\eps)$ lifts vertically as a graph containing $p$, this lift is contained in 
$X_{\eps_0-\eps}$.

A component $c$ of $S_n\cap \partial D(\bar q_j,\eta)\times[-\eps_0,0]$ meeting 
$X_{\eps_0-\eps}$ is a vertical graph, so it is a loop winding around $\bar q_j\times \R$. 
$c$ is called bounding if it bounds in $C_j(\eta,\eps_0)$. We know that if $c$ meets 
$\T\times[-\eps_0,-\eps]$ it bounds a disk so a non bounding curve is contained in 
$X_\eps$. We remark that $S_n$ has a finite number of non bounding curves $c$. 

If there is no non-bounding curve, each boundary curve of $\widetilde S_n$ bounds 
some component. So adding these components to $\widetilde S_n$, we construct a 
connected component of 
$S_n$ in $X_{\eps_0-\eps}$ with no boundary. This is impossible since the mean 
curvature vector of $T_t$ points up.

To any non bounding curve $c$, we can associate a lift of $T_0\setminus \cup_j 
C_j(\eta,\eps)$ containing $c$. Let $F_n$ be the union of the finite number of lifts 
associated to the non bounding curves. Actually if $c$ is a boundary component of $F_n$ 
either it bounds in $C_j(\eta,\eps_0)$, or, if $c$ is non bounding , it is one component of 
the boundary of a connected part of $S_n$ in $C_j(\eta,\eps)$ whose other boundary 
components are 
also in $\partial F_n$. So adding to $F_n$ a finite number of components of $S_n\cap 
\cup_jC_j(\eta,\eps_0)$, we produce a compact minimal surface with no boundary in 
$X_{\eps_0}$. As above such a surface can not exist by the mean convexity of $T_t$.
\end{proof}

%
%
%

\subsection*{Other modifications in the paper}

In~\cite{CoHaMaRo}, Proposition~8 is used in the proof of Theorem~A. Let us explain how 
we have to modify its proof. Let us consider $\boN$ a finite volume hyperbolic 
$3$-manifold. As in the original proof we cut the cusps along constant mean curvature 
tori whose size is sufficiently small such that the above Proposition~\ref{prop} applies 
with $i_0=1$.

As in the original proof, we then glue solid tori with a metric 
depending on a parameter $L$ in order to get a compact manifold $\widetilde \boN$ with 
a metric $g_L$. Moreover the width of $(\widetilde N,g_L)$ is uniformly bounded by $L$. 
Then we can apply Theorem~A in~\cite{MaRo2} to produce a minimal surface $\Sigma$ in 
$(\widetilde{\boN},g_L)$ of index at most $1$ and area at most $M_0$. Then the same 
argument as in the original proof can be applied to prove that, for large $L$, $\Sigma$ 
lies actually in $\boN$. Proposition~\ref{prop} applies since we know that $\Sigma$ has 
index at most $1$.

Actually, this change has also an impact on the proof of Theorem~26 in~\cite{MaRo2}. 
The modification is that in the proof we only have to consider minimal surfaces with 
index at most $1$. So the arguments are the same as those of the previous paragraph.

\bibliographystyle{plain}
\bibliography{../reference.bib}

\end{document}